\documentstyle[franklin]{elsart}

\def\tra{\mathop{\scriptscriptstyle \top}\nolimits} 
\def\smp{\mathop{\scriptscriptstyle +}\nolimits}    
\def\smm{\mathop{\scriptscriptstyle -}\nolimits}
\def\_#1{\mathop{\hspace{-2pt}^{}_{#1}}} 
\def\suml{\mathop{\sum}\limits}
\def\mize{\mathop{\rm minimize}\limits}
\def\conve{\mathop{\longrightarrow}\limits}
\def\cdc{,\ldots,}
\def\on{1,\ldots,n}
\def\om{1,\ldots,m}
\def\s{\mathop{\cal S}\nolimits}
\def\A{\mathop{\cal A}\nolimits}
\def\O {\mathop{\vec{1}}\nolimits}
\def\ww{\mathop{\vec{w}}\nolimits}
\def\ll{\mathop{\vec{l}}\nolimits}
\def\ss{\mathop{\vec{s}}\nolimits}
\def\si{s\_i}
\def\sj{s\_j}
\def\sk{s\_k}
\def\RR{\Rset}
\def\ve{\mathop{\varepsilon}\nolimits}
\def\la{\lambda}
\def\l{\ell}
\newcommand{\ie}{i.e.}
\newcommand{\eg}{e.g.}
\def\SCM{Self-Consis\-tent Monotonicity}
\def\Scm{Self-consistent monotonicity}
\def\scm{self-consistent mono\-to\-ni\-city}
\begin{document}
\def\@date{30 September 1997}
\def\date{30 September 1997}

\begin{frontmatter}
\title{Preference fusion when the number\\ of alternatives
       exceeds two: indirect scoring procedures$\!$\thanksref{thank}}
\thanks[thank]{%
This work was supported by Russian Foundation for Basic Research Grant
No.\ 96--01--01010. Partial research support from the European
Community under Grant ACE--91--R02 is also gratefully acknowledged.
We are thankful to J.~W.~Moon, D.~S.~Schmerling and
J.~L.~Thompson for their help in our search of literature on indirect
scoring procedures.
}

\author{P. Yu. Chebotarev$\!$\thanksref{NSF}} and
\author{Elena Shamis}
\address{Institute of Control Sciences of the Russian Academy of Sciences\\ 65 Profsoyuznaya Str.,
Moscow 117997, Russia%
}
\thanks[NSF]{$\!\!${\tt pchv@rambler.ru}; Participation in the Workshop on Foundations of
Information / Decision Fusion was supported by NSF.}

\begin{abstract}
We consider the problem of aggregation of incomplete preferences
represented by arbitrary binary relations or incomplete paired
comparison matrices. For a number of indirect scoring procedures we
examine whether or not they satisfy the axiom of {\em \scm}. The class
of {\em win-loss combining scoring procedures\/} is introduced which
contains a majority of known scoring procedures. Two main results are
established. According to the first one, every win-loss combining
scoring procedure breaks {\scm}. The second result provides a
sufficient condition of satisfying {\scm}.
\end{abstract}
\end{frontmatter}

\section{Introduction}
\label{Section1}

A good method of preference fusion when there are only two
alternatives is the simple majority of judge votes. This has been
shown by Condorcet in his famous jury theorem \cite{Condorcet1785};
for a review of further developments see
\cite{GrofmanOwen1986,Berg1994}.

The situation becomes much more complicated when the number of
alternatives exceeds two. In this case, the majority of judge
preferences may turn out to be intransitive and contain
preference cycles. Some statistical approaches to this problem
are developed in the theory of paired comparisons (see
\cite{David1988}) and ranking data (see \cite{FlignerVerducci1993}).

In this paper, we use a normative approach typical of the social
choice theory and consider indirect scoring procedures as the methods
of preference fusion. These ingenious procedures are mainly developed
in such disciplines as applied statistics, scoring of sport
tournaments, graph theory, management science, etc.  Indirect score of
an alternative reflects not only the outcomes of its comparisons with
other alternatives, but also the comparison outcomes of those
alternatives to which this one has been compared, and this way it may
depend on all the available preference data. Still there are very few
papers on indirect scoring procedures in the social choice theory.
Probably, the reason is that these procedures are not arithmetically
ascetic enough, and thereby, according to a widespread opinion, do not
correspond to the non-quantitative nature of the social choice
problem.  Another problem is that they are not easy to describe to
the individuals involved in the choice process.  However,
complexity may be considered to be an advantage in the context of the
strategic behavior: complicated procedures are more difficult to
manipulate.

Nevertheless, there is a case where the indirect scoring procedures
are really needed, namely, the case of incomplete preferences. In this
paper, we attempt to demonstrate this case, give a brief critical
review of indirect scoring procedures, and compare their
properties.

Let us consider the situation where an individual is given a set of
alternatives to be compared, but he/she is not an expert in all of
them. Suppose she is allowed to compare only those alternatives she is
familiar with; moreover, to report only those comparison outcomes she
is certain of. In fact, we consider even more general types of
individual opinions, namely, each of them is an arbitrary binary
relation. To prepare these data for calculating indirect scores, we
represent them by {\em incomplete matrices of paired comparisons\/}
borrowed from the statistical theory of paired comparisons.

The main aim of this paper is to examine whether or not the
indirect scoring procedures satisfy the axiom of {\em
Self-Consistent Monotonicity {\rm(}SCM}). It has been introduced
in \cite{ChebotarevShamis1997} where we used it to explore a series of
preference aggregation procedures based on the resolution of discrete
optimization problems.

Up to now there has not been a comprehensive review of indirect
scoring procedures. However, a lot of information can be found in
monographs \cite{David1988,Moon1968,BelkinLevin1990,%
VanBloklandVogelesang1991,CookKress1992} and papers \cite
{Daniels1969,MoonPullman1970,David1971,ShmerlingDubrovskii'1977,%
DavidAndrews1993,Keener1993}.  For relations to the social choice
framework see \cite{LevinNalebuff1995}. An adjacent subject is the
analysis of tournament social choice rules (see, {\eg},
\cite{LaffondLaslier1995,LaffondLaine1996}).

The paper is organized as follows. After the following section which
provides main notation, in Section~\ref{Section3} we give an example
demonstrating the specific character of aggregating incomplete
preferences. Section~\ref{Section4} contains the statement and
justification of SCM. Section~\ref{Section5} demonstrates that if a
procedure is {\em based on individual scores\/}, it break SCM.
Sections~\ref{Section6}, \ref{Section8} and \ref{Section10} give a
review of indirect scoring procedures known from the literature.
Section~\ref{Section7} introduces the class of {\em win-loss combining
scoring procedures\/} and establishes that these procedures fail to
satisfy SCM.  Section~\ref{Section9} provides a sufficient condition
for SCM, and Section~\ref{Section11} discusses two other axioms.

\section{Notation}
\label{Section2}

Suppose $J=\{\on\}$ is a set of alternatives to be compared. There are
$m$ individuals (judges, voters, etc.), and each of them reports
his/her preferences. Concerning any pair $(i,j)$ of alternatives, an
individual may give one of the following four responses: (a) ``$i$ is
better than $j$'', (b) ``$j$ is better than $i$", (c) ``$i$ and $j$
are equivalent", and (d) ``I do not report my opinion on this pair."
Using the responses of the $p$th individual, an {\em incomplete paired
comparison matrix\/} $A^{(p)}=[a_{ij}^p]$ is filled out. Its entries
$a_{ij}^p$ and $a_{ji}^p$ correspond to the comparison outcome of $i$
and $j$ as follows:

if the $p$th individual said that
\begin{itemize}
\item[--] $i$ is better than $j$, then $a_{ij}^p=1$, $a_{ji}^p=0$;
\item[--] $i$ and $j$ are equivalent, then $a_{ij}^p=a_{ji}^p=1/2$;
\end{itemize}
If the individual did not report his/her opinion, $a_{ij}^p$ and
$a_{ji}^p$ remain undefined.

By definition, we set all diagonal entries to be zero:  $a_{ii}^p=0$,
$i=\on$, $p=\om.$ Some authors put $a_{ii}^p=1/2$ or leave $a_{ii}^p$
undefined, but this does not change the results significantly.

This way we have an array of incomplete paired comparison matrices
$\A=(A^{(1)}\cdc A^{(m)})$ which is referred to as a {\em profile of
individual preferences}. A {\em scoring\/} ({\em rating}) {\em
procedure\/} is a function $\s$ from the set of all profiles (or its
subset) into $\RR^n$, where $s\_i$, the $i$th component of the
resulting column vector $\ss=(s\_1\cdc s\_n)^{\tra}$ is interpreted as
a {\em score of the alternative $i$}.  All the scoring procedures
considered in this paper are assumed to be {\em neutral\/} (any
reindexing of the alternatives preserves their scores) and {\em
anonymous\/} (any reindexing of the individuals preserves the scores of
the alternatives). The term ``aggregation of incomplete preferences"
will mean here nothing but rating the alternatives by the scores
$s\_1\cdc s\_n$: the greater is a score, the more socially preferred is
the alternative.

\section{An Example}
\label{Section3}

The problem of aggregation of incomplete preferences has a specific
flavor which cannot be revealed in dealing with complete preferences.
To capture it consider the following example.

\begin{figure}[htb]
\setlength{\unitlength}{5mm}
\begin{picture}(27.2,8)
\put(0, 1){\line(1,0){6}}
\put(0, 7){\line(1,0){6}}
\put(0, 1){\line(0,1){6}}
\put(6, 1){\line(0,1){6}}
\put(0, 1){\begin{picture}(6,6)
             \put(1.4,1.4){2}
             \put(1.4,4.4){1}
             \put(4.4,2.0){4}
             \put(4.4,3.8){3}
             \put(1.3,4.0){\vector(0,-1){1.9}} 
             \put(1.6,4.0){\vector(0,-1){1.9}} 
             \put(1.9,4.0){\vector(0,-1){1.9}} 
             \put(2.8,-1){\small a.}
          \end{picture}}
\put( 7, 1){\line(1,0){6}}
\put( 7, 7){\line(1,0){6}}
\put( 7, 1){\line(0,1){6}}
\put(13, 1){\line(0,1){6}}
\put( 7, 1){\begin{picture}(6,6)
             \put(1.4,1.4){2}
             \put(1.4,4.4){1}
             \put(4.4,2.0){4}
             \put(4.4,3.8){3}
             \put(2.0,4.3){\line(4,-1){2.0}} 
             \put(2.0,4.6){\line(4,-1){2.0}} 
             \put(2.0,4.9){\line(4,-1){2.0}} 
             \put(2.8,-1){\small b.}
          \end{picture}}
\put(14, 1){\line(1,0){6}}
\put(14, 7){\line(1,0){6}}
\put(14, 1){\line(0,1){6}}
\put(20, 1){\line(0,1){6}}
\put(14, 1){\begin{picture}(6,6)
             \put(1.4,1.4){2}
             \put(1.4,4.4){1}
             \put(4.4,2.0){4}
             \put(4.4,3.8){3}
             \put(2.0,1.3){\line(4,1){2.0}} 
             \put(2.0,1.6){\line(4,1){2.0}} 
             \put(2.0,1.9){\line(4,1){2.0}} 
             \put(2.8,-1){\small c.}
          \end{picture}}
\put(21, 1){\line(1,0){6}}
\put(21, 7){\line(1,0){6}}
\put(21, 1){\line(0,1){6}}
\put(27, 1){\line(0,1){6}}
\put(21, 1){\begin{picture}(6,6)
             \put(1.4,1.4){2}
             \put(1.4,4.4){1}
             \put(4.4,2.0){4}
             \put(4.4,3.8){3}
             \put(1.3,4.0){\vector(0,-1){1.9}} 
             \put(1.6,4.0){\vector(0,-1){1.9}} 
             \put(1.9,4.0){\vector(0,-1){1.9}} 
             \put(2.0,1.3){\line(4,1){2.0}} 
             \put(2.0,1.6){\line(4,1){2.0}} 
             \put(2.0,1.9){\line(4,1){2.0}} 
             \put(2.0,4.3){\line(4,-1){2.0}} 
             \put(2.0,4.6){\line(4,-1){2.0}} 
             \put(2.0,4.9){\line(4,-1){2.0}} 
             \put(2.8,-1){\small d.}
          \end{picture}}
\end{picture}
\caption{a, b, c: Preferences of nine individuals; d: Combined
preferences.}
\label{Fig1}
\end{figure}
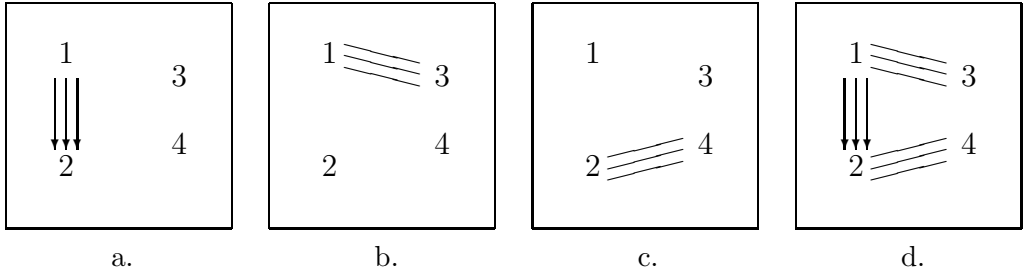

Suppose there are four alternatives, $J=\{1,2,3,4 \}$, and nine
individuals (judges).  Three judges are familiar with alternatives 1
and 2 and they all feel that 1 is better (Fig.~\ref{Fig1}.a where an
arrow from 1 to 2 shows that 1 is preferred to 2). Three judges think
that 1 and 3 are equivalent (Fig.~\ref{Fig1}.b where line segments
without arrows designate equivalencies).  Finally, three judges are
familiar with 2 and 4, and they all consider these alternatives of
equal quality (Fig.~\ref{Fig1}.c).  The preferences of these nine
individuals are incorporated in Fig.~\ref{Fig1}.d.  Can we say
anything about the comparative quality of 3 and 4? Although no
judge is familiar with 3 and 4 simultaneously, it is plausible that 3
is better than 4.  Note that 3 and 4 have the same comparison
outcomes:  three equivalencies; they differ only in the alternatives
to which they were compared. The ``opponent" of 3 is probably better
than that of 4 and due to this 3 may be estimated higher than 4.

Thus, in aggregating incomplete preferences we should take into account
``quality (strength) of the opponents" or ``caliber of the opposition"
or ``schedule difficulty". These terms came from sport, and we see that
there is a similarity between the problem of aggregation of
incomplete preferences and the problem of rating the participants of an
incomplete tournament (note that the extent of this similarity is a
distinct and interesting question).  Moreover, in the following section
this sport analogy will be exploited to explain the meaning of
Self-Consistent Monotonicity, the axiom which is discussed throughout
the paper.  Sometimes the comparison outcome where $i$ is preferred to
$j$ will be called a ``{\em win}" of $i$ and a ``{\em loss}" of $j$.

\section{{\SCM}}
\label{Section4}

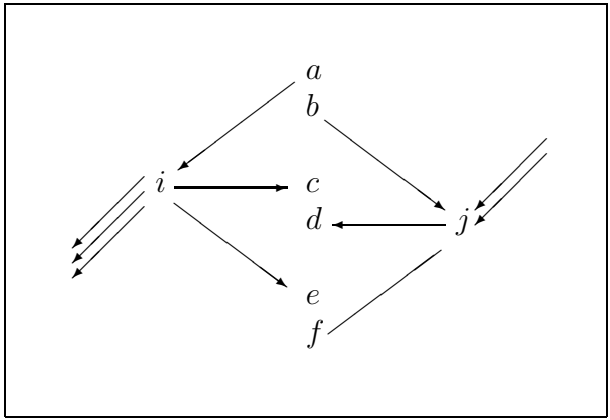
\begin{figure}[htb]
\setlength{\unitlength}{5mm}
\begin{picture}(25,12)
\put( 5, 0){\line(1,0){16}}
\put( 5,11){\line(1,0){16}}
\put( 5, 0){\line(0,1){11}}
\put(21, 0){\line(0,1){11}}
\put( 5, 0){\begin{picture}(16,11)
             \put(4.0,6.0){$i$}
             \put(8.0,2.0){$f$}
             \put(8.0,3.0){$e$}
             \put(8.0,5.0){$d$}
             \put(8.0,6.0){$c$}
             \put(8.0,8.0){$b$}
             \put(8.0,9.0){$a$}
             \put(12.0,5.0){$j$}
             \put(7.7,8.9){\vector(-4,-3){3.1}} 
             \put(4.5,6.1){\vector(1  ,0){3.0}} 
             \put(4.5,5.7){\vector(4, -3){3.0}} 
             \put(8.5,7.9){\vector(4 ,-3){3.2}} 
             \put(11.7,5.1){\vector(-1 ,0){3.0}} 
             \put(8.6,2.2){\line  (4,  3){3.0}} 
             \put(3.7,5.6){\vector(-1,-1){1.9}} 
             \put(3.7,6.0){\vector(-1,-1){1.9}} 
             \put(3.7,6.4){\vector(-1,-1){1.9}} 
             \put(14.4,7.0){\vector(-1,-1){1.9}} 
             \put(14.4,7.4){\vector(-1,-1){1.9}} 
          \end{picture}}
\end{picture}
\caption{An illustration to {\scm}.}
\label{Fig2}
\end{figure}

First, let us consider one more example. Suppose we have the incomplete
sport tournament whose fragment containing all game results of $i$ and
$j$ is shown in Fig.~\ref{Fig2}.  Suppose we know that $a$ is stronger
than $b$, $c$ is stronger than $d$, and $e$ is stronger than $f$. Then
the results of $i$ are definitely better than those of $j$. Indeed, a
loss to $a$ is more pardonable than that to $b$, a win over $c$ is
more honorable than that over $d$, and a win over $e$ is more valuable
than a draw with $f$.  Besides, $i$ has three extra wins and $j$ has
two extra losses, which intensify the advantage of $i$ over $j$.

In fact, the following {\scm} axiom requires that in such situations
the score of $i$ should be greater than the score of $j$. The only
point is the meaning of preconceptions like ``we know that $a$ is
stronger than $b$". {\Scm} applies to scoring procedures, and here ``we
know that $a$ is stronger than $b$" signifies ``{\em this\/} scoring
procedure gives $a$ a greater score than $b$." This makes clear why
this kind of monotonicity is called ``self-consistent".

In the following statement of {\scm}, we use the term ``multiset". Its
difference from set is that multiset may contain the same element in
several copies. This concept is needed here since the outcomes of
comparisons with different alternatives may coincide. The elements of
multisets will be written within angle brackets.

\medskip
{\bf {\Scm} (SCM).} {\em
A scoring procedure $\s$ is {\em Self-Consistently Monotonic\/} if for
any set of alternatives $J$, any profile $\A$, and any alternatives
$i,j\in J$ it satisfies the following condition.
\par
Suppose $U_i=\langle a_{ik}^p\mid k,p\rangle$ and $U_j=\langle
a_{j\ell}^q\mid \ell,q\rangle$ are the multisets of the
comparison outcomes of alternatives $i$ and $j$, respectively.
Suppose that $U_i$ can be split into $U_i^I$ and $U_i^{II}$, and
$U_j$ can be split into $U_j^I$ and $U_j^{II}$ in such a way that
\begin{itemize}
\item[{\rm(}i{\rm)}] $a_{ik}^p\in U_i^I$ implies $a_{ik}^p=1$,
      $a_{j\l}^q\in U_j^I$ implies $a_{j\l}^q=0$;
\item[{\rm(}ii{\rm)}] there exists a one-to-one correspondence $\pi$
      of $U_i^{II}$ onto $U_j^{II}$ such that
      $\pi(a_{ik}^p)=a_{j\l}^q$ implies $a_{ik}^p\ge a_{j\l}^q$ and
      $\sk\ge s\_\l$.
\end{itemize}
\par
Then $\si\ge \sj$.
\par
Moreover, if in addition $U_i^I$ is nonempty or $U_j^I$ is nonempty or
at least one inequality in (ii) is strict at least once, then
$\si>\sj$.
}
\bigskip

The application of SCM to a pair of alternatives $(i,j)$ will be called
{\em confrontation of $i$ and $j$}.

It turns out that many well-known indirect scoring procedures fail to
satisfy {\scm}. In the following two sections we try to explore the
features of these procedures that cause them not to satisfy SCM. The
results are presented in Proposition~\ref{Prop1}
(Section~\ref{Section5}) and Theorem~\ref{Theorem1}
(Section~\ref{Section7}). Section~\ref{Section9} provides a sufficient
condition of SCM (Theorem~\ref{Theorem2}).

\section{Procedures Based on Individual Scores}
\label{Section5}

A scoring procedure $\s$ is {\em based on individual scores\/} if there
exist functions $f$ and $\delta$ such that for any profile
$\A=(A^{(1)}\cdc  A^{(m)})$, the corresponding score vector $\ss$ can
be expressed as $\ss=\delta(\ss^{(1)}\cdc  \ss^{(m)})$, where
$\ss^{(p)}$ is a partial score vector depending solely on the
comparison matrix of individual $p$: $\ss^{(p)}=f(A^{(p)})$, $p=\om$.

The most important instance is
$$
\ss=\suml_{p=1}^m \ss^{(p)},
$$
forming (in case of complete preferences) a class of procedures which
satisfy the recent axiomatics by Myerson \cite{Myerson1995}.

\begin{prop}\label{Prop1}
There are scoring procedures based on individual scores that satisfy
SCM on complete preferences but no such procedure satisfies SCM for
incomplete preferences.
\end{prop}

\begin{pf}
To prove the first statement, it suffices to consider the row sum
procedure, which in case of complete preferences has the form
\begin{eqnarray}
s_i^{(p)}&=&\suml_{j=1}^n a_{ij}^p,\qquad i=\on,\quad p=\om,
\nonumber\\
\ss&=&\suml_{p=1}^m \ss^{(p)}.
\label{case}
\nonumber
\end{eqnarray}
This procedure obviously satisfies SCM.

The second statement can be proved by considering our first example of
Fig.~\ref{Fig1}. Let $p=1,2,3$ be the numbers of individuals in
Fig.~\ref{Fig1}.a, $p=4,5,6$ in Fig.~\ref{Fig1}.b, and $p=7,8,9$ in
Fig.~\ref{Fig1}.c. By neutrality and anonymity of $\s$,

\begin{eqnarray*}
&s_3^{(1)}=s_3^{(2)}=s_3^{(3)}=s_4^{(1)}=s_4^{(2)}=s_4^{(3)},&\\
&s_3^{(4)}=s_3^{(5)}=s_3^{(6)}=s_4^{(7)}=s_4^{(8)}=s_4^{(9)},&\\
&s_4^{(4)}=s_4^{(5)}=s_4^{(6)}=s_3^{(7)}=s_3^{(8)}=s_3^{(9)}.&
\end{eqnarray*}
Hence, $s\_3=\suml_{p=1}^9 s_3^{(p)}=\suml_{p=1}^9 s_4^{(p)}=s\_4.$

Confronting 3 and 4 by SCM (see Fig.~\ref{Fig1}.d), we deduce
$s\_1=s\_2.$ But now, confronting 1 and 2, we get contradiction.
\qed
\end{pf}

In words, these procedures generally break SCM because the result of
confrontation between two alternatives may be only determined by the
preferences of the {\em whole\/} board of judges.

\section{Aggregate Scoring Procedures}
\label{Section6}

Hereafter we will consider scoring procedures that cannot be reduced to
individual scores. Such indivisible procedures can be called {\em
aggregate}. A large variety of them are representable through the
matrix $A=[a\_{ij}]$ of total comparison outcomes on the pairs of
alternatives:
\begin{equation}
a\_{ij}=\cases{0, &if $a_{ij}^p$ is undefined for every $p,$\cr
               \suml_{p}a_{ij}^p, & otherwise.\cr}
\label{0}
\end{equation}

Consider the scores of the form:
\begin{equation}
s\_i=\suml_{j=1}^n f(a\_{ij},x\_j), \quad i=\on,
\label{1}
\end{equation}
where $x\_j$ $(j=\on)$ is some estimate of the alternative $j$, $f$ is
a function nondecreasing in both $a\_{ij}$ and $x\_j$. Sometimes
averaged scores of the form
\begin{equation}
s\_i={1\over m\_i}\suml_{j=1}^n f(a\_{ij},x\_j), \quad i=\on
\label{1.1}
\end{equation}
are proposed, where $m\_i$ is the total number of comparisons of $i$
(see, {\eg}, \cite{Keener1993}). For simplicity, we do not consider
this modification here, however, it is covered by
Theorem~\ref{Theorem1} below.

\subsection{Wei's Procedure}
\label{Subsect_Wei}

The most popular form of $f$ in Eqs~(\ref{1}) is the product:
\begin{equation}
s\_i=\suml_{j=1}^n a\_{ij}x\_j, \quad i=\on.
\label{2}
\end{equation}

Regarding the choice of the form of $x\_j$, perhaps the most attractive
idea is to relate $x\_j$ and $s\_j$ directly. However, if we set
$x\_j=\sj,$ $j=\on$, the homogeneous system of linear equations
(\ref{2}) may have only a trivial solution $s\_1=\cdots=s\_n=0.$ A
minor modification is to define $x\_j$ to be proportional to $\sj$:
$$
x\_j=\sj/\la, \quad j=\on,
$$
not specifying $\la$ {\it a priori}. Then
\begin{equation}
\la s\_i=\suml_{j=1}^n a\_{ij}s\_j, \quad i=\on,
\label{3}
\end{equation}
or in the matrix form,
\begin{equation}
\la\ss=A\ss,
\label{4}
\end{equation}
thus $\la$ is an eigenvalue and $\ss$ an eigenvector of $A$.

This scoring procedure was proposed by Wei \cite{Wei1952} and
became well-known after Kendall's \cite{Kendall1955} paper. Matrix
$A$ generally has several eigenvalues, each having its own subspace of
eigenvectors.  Only one solution is taken, and this solution possesses
special properties.  Let us suppose that the preferences are {\em
indivisible}, {\ie}, the set of alternatives $J$ cannot be split into
two nonempty parts $J_1$ and $J_2$ such that for no alternatives $j\in
J_2$ and $i\in J_1,$ $a\_{ji}>0.$ Then, by the Perron--Frobenius
theorem, the largest in absolute value eigenvalue of $A$ is positive,
the corresponding subspace of eigenvectors is one-dimensional, and
these eigenvectors have the same sign of all terms. Just the normalized
positive eigenvector $\ss$ of this type is taken as a vector of scores
in the Wei procedure (which applies to indivisible preferences).

Wei's eigenvector can be obtained iteratively. Consider the following
sequence of score vectors:
\begin{eqnarray*}
\ss^1&=&A\O, \quad {\rm where}\;\:\O=(1\cdc 1)^{\tra};\\
\ss^2&=&A\ss^1=A^2\O;\\
& &\cdots\cdots\\
\ss^k&=&A\ss^{k-1}=A^k\O;\\
& &\cdots\cdots
\end{eqnarray*}

It follows from the Perron--Frobenius theorem that in case of
indivisible profiles, the normalized sequence $(\ss^k)$ converges to
Wei's score vector:
\begin{equation}
{\ss^k\over \suml_{i=1}^n s_i^k}={1\over \suml_{i=1}^n s_i^k}
A^k\O\conve_{k\to\infty}\ss.
\label{5}
\end{equation}

It is a simple fact that $s_i^k$ is the number of $k$-length paths
from $i$ to all vertices in the preference multidigraph. Thus, in
Wei's procedure, alternatives are compared by the number of very long
(``infinitely-long") paths diverging from them in the preference
multidigraph. In Subsection~\ref{Subsect_Katz} we shall consider a
scoring procedure where longer paths are accounted with smaller
weights.

\begin{prop}\label{Prop2}
\newcounter{Prop2}\setcounter{Prop2}{\value{thm}}
If indivisible preferences are complete, Wei's scoring procedure
satisfies SCM, but it breaks SCM  for incomplete indivisible
preferences.
\end{prop}

\begin{pf}
To prove the first statement, confront arbitrary $i$ and $j$ by SCM:
\begin{equation}
s\_i={1\over\la}\suml_{k=1}^n a\_{ik}\sk,
\label{6}
\end{equation}
\begin{equation}
s\_j={1\over\la}\suml_{k=1}^n a\_{jk}\sk,
\label{7}
\end{equation}
and assume that some partitions of $U_i$ and $U_j$ described in {\scm}
exist. Then $U_i^{II}$ induces on the right-hand side of Eq.~(\ref{6})
a no lesser (respectively, greater in the strict case) sum of terms
than $U_j^{II}$ induces on the right-hand side of Eq.~(\ref{7}).  The
same can be said of $U_i^{I}$ and $U_j^{I}$.  Therefore, $\si\ge\sj$
($\si>\sj$ in the strict case).

To prove the second statement, consider the cumulative preference
digraph depicted in Fig.~\ref{Fig3}.

\begin{figure}[htb]
\setlength{\unitlength}{5mm}
\begin{picture}(19,8.6)
\put(9, 0){\line(1,0){10}}
\put(9, 8){\line(1,0){10}}
\put(9, 0){\line(0,1){8}}
\put(19,0){\line(0,1){8}}
\put(9, 0){\begin{picture}(10,8)
             \put(2.0,6.0){$1$}
             \put(2.9,3.8){$3$}
             \put(3.9,1.9){$4$}
             \put(6.0,5.0){$2$}
             \put(8.7,3.6){$5$}
             \put(2.3,5.7){\vector( 1,-2){0.6}} 
             \put(3.3,3.7){\vector( 1,-2){0.6}} 
             \put(5.9,4.8){\vector(-2,-3){1.53}} 
             \thicklines
             \put(7.3,3.9){\line  (1,  0){1.2}} 
             \thinlines
            \put(4.0,4.0){\oval(6,6)}
          \end{picture}}
\end{picture}
\caption{Preferences in the proof of Proposition~$\!$\arabic{Prop2}.}
\label{Fig3}
\end{figure}

The comparison outcomes of alternative 5 are equivalencies with
alternatives 1, 2, 3 and 4. Alternative 5 is introduced only to make
the preferences indivisible. Consider the equations with $s\_2$ and
$s\_3$ on the left-hand sides:
\begin{eqnarray*}
\la s\_2&=&s\_4+s\_5/2,\\
\la s\_3&=&s\_4+s\_5/2.
\end{eqnarray*}
Hence, $s\_2=s\_3$, which contradicts SCM, since 3 has an extra ``loss"
($a\_{31}=0,\;{a\_{13}=1}$).
\qed
\end{pf}

The example used in the proof of Proposition~\ref{Prop2} reveals one
important feature of Wei's procedure as applied to incomplete
preferences. Namely, this procedure bases itself only on ``wins" and
does not take into account ``losses". More precisely, it does not
distinguish ``losses" from missing comparisons (and precisely this
contradicts SCM). A possible way to make the procedure more balanced is
to also consider a dual procedure which is based on ``losses", and
then to combine the resulting scores. This has been done by Hasse
\cite{Hasse1961} and Ramanujacharyulu \cite{Ramanujacharyulu1964}.

\subsection{The Procedures by Hasse and Ramanujacharyulu}
\label{Subsect_Has}

Let $\ww$ and $\ll$ denote the right and left eigenvectors of matrix
$A$, corresponding to its maximal eigenvalue $\la$ :
\begin{equation}
\la w\_i=\suml_{j=1}^n a\_{ij}w\_j, \quad i=\on,
\label{8}
\end{equation}
\begin{equation}
\la \l\_i=\suml_{j=1}^n a\_{ji}\l\_j, \quad i=\on,
\label{9}
\end{equation}
Here $\ww=(w\_1\cdc w\_n)$ is the vector of Wei's {\em win-scores},
which is not sensitive to the distinction between ``losses" and missing
comparisons; $\ll=(\l\_1\cdc \l\_n)$ is the vector of {\em
loss-scores}, which does not distinguish ``wins" and missing
comparisons. The matrix form of Eqs~(\ref{9}) is
\begin{equation}
\la\ll=B\ll,
\label{10}
\end{equation}
where $B=A^{\tra}$ is the transpose of $A$.

Hasse \cite{Hasse1961} proposed to combine $\ww$ and $\ll$ into the
ultimate scores:
\begin{equation}
\si=w\_i-\l\_i,\quad i=\on.
\label{11}
\end{equation}
Another procedure proposed by Ramanujacharyulu
\cite{Ramanujacharyulu1964} has scores given by:
\begin{equation}
\si=w\_i/\l\_i,\quad i=\on.
\label{12}
\end{equation}

These two procedures are close but may produce different rankings,
since the first one registers the absolute differences between
win-scores and loss-scores, whereas the second is based on the relative
differences. Both procedures are applicable to indivisible preferences
and turn out to break {\scm}. We do not formulate the corresponding
statement, because there are many possible modifications of such
procedures, and a more reasonable thing to do is to isolate their
common features and then to prove a more general theorem.  To this
end, in the following two subsections we consider two other families
of indirect scoring procedures, and then proceed with a theorem.

\subsection{The Procedure of Katz--Thompson--Taylor}
\label{Subsect_Katz}

As we have seen, Wei's scores rank order the alternatives according
to the number of very long paths diverging from them in the preference
multidigraph. Katz \cite{Katz1953} and then Thompson
\cite{Thompson1958} and Taylor \cite{Taylor1968} proposed a scoring
procedure in which the long paths play just a correcting role, whereas
the short paths are taken into account with greater weights. Namely,
Thompson introduced the scores of the form
\begin{equation}
\ww=(A+\ve A^2+\ve^2 A^3+\cdots)\O,
\label{13}
\end{equation}
where $\ve$ is a small positive parameter. If $\ve<r^{-1}$ where $r$
is the spectral radius of $A$, then the series (\ref{13}), which in
fact is a geometric progression, converges and
\begin{equation}
\ww=A(I-\ve A)^{-1}\O,
\label{14}
\end{equation}
where $I$ is the identity matrix.

Thompson shows that his normalized vector $\ww$ tends to Wei's
eigenvector as $\ve$ approaches $r^{-1}$ from below. Besides, he gives
a game-theoretical interpretation of these scores comparable with the
approach by Laffond, Laslier and Lebreton \cite{LaffondLaslier1993}.

The paper by Taylor is one of the few social choice articles devoted
to indirect scoring procedures. He investigated sequential voting by
the committee and considered  Hamiltonian sequences of winning
alternatives, which led him to a scoring procedure very close to
(\ref{13}). In fact, Taylor used another matrix $C=[c\_{ij}]$ instead
$A$, where
$$
c\_{ij}=\cases{a\_{ij}-a\_{ji}, &if $a\_{ij}-a\_{ji}\ge0$,\cr
                             0, & otherwise.\cr}
$$
Since for the preference structure we use in the proof of
Theorem~\ref{Theorem1} in Section~\ref{Section7}, $A$ and $C$
coincide, the statement of this theorem is valid for Taylor's
procedure.

As above, $\ww$ in Eqs~(\ref{13}) and (\ref{14}) is a vector of {\em
win-scores}, and it is worth being supplemented by the corresponding
{\em loss-scores}. So we may consider a pair of vectors $(\ww,\ll)$
where
\begin{equation}
\ll=(B+\ve B^2+\ve^2 B^3+\cdots)\O,
\label{14a}
\end{equation}
and combine $\ww$ and $\ll$ as in the procedures of Hasse and
Ramanujacharyulu:
$$
\si=w\_i-\l\_i,\quad i=\on
$$
or
$$
\si=w\_i/\l\_i,\quad i=\on.
$$

The question we are interested in is whether or not these procedures
satisfy {\scm}. The answer is negative. After considering one more
family of scoring procedures, we give a theorem that will clarify
this point to some extent.

\subsection{Directed Tree Procedure by
Daniels--Ushakov--Goddard--Levchenkov}
\label{Subsect_Daniels}

Let $c_i^{\smp}$ and $c_i^{\smm}$ be the total ``win" and total
``loss" of $i$, respectively:
\begin{eqnarray}
c_i^{\smp}&=&\suml_{j=1}^n a\_{ij},
\nonumber\\
c_i^{\smm}&=&\suml_{j=1}^n a\_{ji},\quad i=\on.
\label{15}
\end{eqnarray}
Daniels \cite{Daniels1969} and other writers \cite{Ushakov1971,%
Ushakov1976,Goddard1983,Levchenkov1990,Levchenkov1992} proposed, with
various motivations, the scores $w\_1\cdc w\_n$ that satisfy
\begin{equation}
w\_i={1\over c_i^{\smm}}\suml_{j=1}^n a\_{ij}w\_j,\quad i=\on.
\label{16}
\end{equation}
The only difference of this system of equation from that of Wei is
that $\la$ is replaced by $c_i^{\smm}$. This is possible because the
matrix $[a\_{ij}/c_i^{\smm}]$ always has a unit eigenvalue. Scores
$w\_i$ have several interesting interpretations listed below.
\begin{itemize}
\item[(i)] $w\_i$ is proportional to the number of directed trees
diverging from $i$ in the preference multidigraph \cite{Berman1980};
\item[(ii)] $w\_i$ are proportional to the final probabilities of the
``winning Markov chain" \cite{Ushakov1971,Ushakov1976,Daley1979,%
Levchenkov1990,Levchenkov1992,Levchenkov1994};
\item[(iii)] $\ww$ is the vector of ``fair bets" in the model of
Moon and Pullman \cite{MoonPullman1970}.
\end{itemize}

Also, in \cite{LevchenkovGrivko1992,Levchenkov1993a,%
GrivkoLevchenkov1994} this procedure is analyzed in the social choice
framework, and in \cite{Levchenkov1993b,Levchenkov1995} its strategic
properties are touched upon.

Exactly as before, we may consider the corresponding
{\em loss-scores\/} mentioned by Daley \cite{Daley1979}:
\begin{equation}
\l\_i={1\over c_i^{\smp}}\suml_{j=1}^n a\_{ji}\l\_j,\quad i=\on
\label{17}
\end{equation}
and define the ultimate score vector as a combination of $\ww$ and
$\ll$, {\eg},
$$
\si=w\_i-\l\_i,\quad i=\on
$$
or
$$
\si=w\_i/\l\_i,\quad i=\on.
$$
The procedures of this family break SCM, and now we are in a position
to specify a more general class of indirect scoring procedures that do
not obey SCM. All the above procedures are included.

\section{A General Class of Procedures That Fail to Satisfy \SCM}
\label{Section7}

In this section we define {\em win-loss indices\/} and {\em win-loss
combining procedures\/} and prove that such procedures break {\scm}.

In the following Definitions~\ref{Def1}--\ref{Def3}, profile $\A$ is
fixed; ``vector" means $n$-component column vector.

\begin{defn}\label{Def1}
{\rm(}recursive{\rm)}. {\em Win-loss index of finite order.}
\begin{itemize}
\item[$1^\circ.$] A pair of vectors $(\ww,\ll)$ is a {\em win-loss
      index of zero order\/} if $ \ww=\ll= (1\cdc 1)^{\tra}.$
\item[$2^\circ.$] A pair of vectors $(\ww,\ll)$ is a {\em win-loss
      index of order $k$\/} if there exists a win-loss index
      $(\ww',\ll')$ of order $k-1$ and functions $f$, $g$, $\alpha$ and
      $\beta$ such that
      \begin{eqnarray}
      w\_i&=&{1\over \alpha(c_i^{\smp},c_i^{\smm})}f(W'_i),
      \nonumber\\
      \l\_i&=&{1\over \beta(c_i^{\smm},c_i^{\smp})}g(L'_i),\quad i=\on,
      \label{18}
      \end{eqnarray}
      where $W'_i=\langle(a\_{i1},w'\_1),(a\_{i2},w'\_2)\cdc
      (a\_{in},w'\_n)\rangle$,
      $L'_i=\langle(a\_{1i},\l'\_1),(a\_{2i},\l'\_2),$
      $\ldots, (a\_{ni},\l'\_n) \rangle$, $i=\on$; angle brackets
      designate multisets.
\end{itemize}
$W'_i$ and $L'_i$ are respectively called the {\em win-performance\/}
and the {\em loss-per\-for\-mance\/} of $i$ corresponding to
$(\ww',\ll')$.
\end{defn}

One can see that in the Directed tree procedure of
Subsection~\ref{Subsect_Daniels}, $\alpha(c_i^{\smp},c_i^{\smm})$ is
$c_i^{\smm}$; $\beta(c_i^{\smm},c_i^{\smp})$ is $c_i^{\smp}$; $f$ and
$g$ are sums of products. An essential difference is that in the
Directed tree procedure, $w\_i$  and $\l\_i$ are related with the
win-performance and loss-performance corresponding to the same pair of
vectors $(\ww,\ll)$, not to a win-loss index of the previous order.
Thus, Eqs~(\ref{16}) and (\ref{17}) do not fit Definition~\ref{Def1}.
That pair of vectors $(\ww,\ll)$ satisfies the following definition.

\begin{defn}\label{Def2}
{\em Win-loss index of infinite order.}
A pair
of vectors
$(\ww,\ll)=((w\_1\cdc w\_n)^{\tra},(\l\_1\cdc \l\_n)^{\tra})$ is
a {\em win-loss index of infinite order\/} if there exist
functions $f$, $g$, $\alpha$ and $\beta$ such that
\begin{eqnarray}
w\_i & = & {1\over \alpha(c_i^{\smp},c_i^{\smm})}f(W_i),
\nonumber\\
\l\_i & = & {1\over \beta(c_i^{\smm},c_i^{\smp})}g(L_i),\quad i=\on,
\label{19}
\end{eqnarray}
where
$
W_i\!=\!\left<(a\_{i1},w\_1),(a\_{i2},w\_2)\cdc (a\_{in},w\_n)\right>
$
and
$
L_i=\langle(a\_{1i},\l\_1),(a\_{2i},\l\_2)\cdc $
$(a\_{ni},\l\_n)\rangle
$,
$i=\on.$
\end{defn}

The following definition introduces win-loss indices without order.

\begin{defn}\label{Def3}
{\em Win-loss index.}
A pair
of vectors
$(\ww,\ll)=((w\_1\cdc w\_n)^{\tra},$ $(\l\_1\cdc \l\_n)^{\tra})$
is a {\em win-loss index\/}
if there exists a sequence $((w^1,\l^1),(w^2,\l^2),...)$ of win-loss
indices of finite or infinite order and functions $\varphi$ and $\psi$
such that for $i=\on$,
\begin{eqnarray}
w\_i &=& \varphi(w_i^1,\l_i^1,w_i^2,\l_i^2,\ldots),
\nonumber\\
\l\_i &=& \psi(w_i^1,\l_i^1,w_i^2,\l_i^2,\ldots).
\label{lemon}
\end{eqnarray}
\end{defn}

The components of $\ww$ and $\ll$ are referred to as {\em
win-scores\/} and {\em loss-scores}, respectively (note that they
depend on both $w$-components and $\l$-components of the sequence
elements). The next (and last) definition introduces a scoring
procedure that combines win-scores and loss-scores of a win-loss index.

\begin{defn}\label{Def4}
{\em Win-loss combining scoring procedure}.
Scoring procedure $\s$ is a {\em win-loss combining scoring
procedure\/} if for any profile $\A$ there exists a win-loss index
$(\ww,\ll)$ and a function $h$ such that
$$
s\_i=h(w\_i,\l\_i), \quad i=\on.
$$
\end{defn}

\begin{prop}\label{Prop3}
All scoring procedures that generate score vectors of
Section~\ref{Section6} or any their iterative approximations
constructed as in $2^\circ$ of Definition~\ref{Def1}, are win-loss
combining scoring procedures.
\end{prop}

The proof of Proposition~\ref{Prop3} is straightforward. Among other
win-loss combining scoring procedures, let us mention the elaborate
procedure by David (see \cite{David1987,David1988,DavidAndrews1993},
and procedures by Cowden \cite{Cowden1975} and Ginovker
\cite{Ginovker1981}.  The latter two authors also proposed other
procedures which will be mentioned below.  An interesting discussion of
related topics can be found in \cite{Keener1993}.

\begin{thm}\label{Theorem1}
\newcounter{Theorem1}\setcounter{Theorem1}{\value{thm}}
Every win-loss combining scoring procedure defined on the set of
indivisible preference profiles breaks {\scm}.
\end{thm}

\begin{pf}
Assume that there exists a win-loss combining scoring procedure that
satisfies SCM on the set of indivisible preference profiles. Suppose
that $(\ww,\ll)$ is its win-loss index for the profile with the
preference digraph shown in Fig.~\ref{Fig4} and that
$s\_i=h(w\_i,\l\_i)$, $i=\on$.

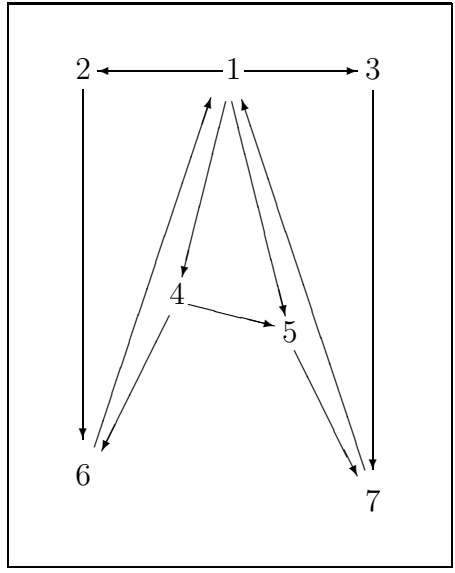
\begin{figure}[htb]
\setlength{\unitlength}{5mm}
\begin{picture}(12,17)
\put(8,1){\begin{picture}(12,13)
     \put(   0, 0){\line(1,0){11.8}}
     \put(   0,15){\line(1,0){11.8}}
     \put(   0, 0){\line(0,1){15}}
     \put(11.8, 0){\line(0,1){15}}
          \put(5.8,13.0){$1$}
          \put(1.8,13.0){$2$}
          \put(9.5,13.0){$3$}
          \put(4.3, 7.0){$4$}
          \put(7.3, 6.0){$5$}
          \put(1.8, 2.2){$6$}
          \put(9.5, 1.5){$7$}
          \put(5.7, 13.2){\vector(-1, 0){ 3.3 }} 
          \put(6.3, 13.2){\vector( 1, 0){ 3.0 }} 
          \put(5.8, 12.4){\vector(-1,-4){ 1.17}} 
          \put(5.95,12.4){\vector( 1,-4){ 1.43}} 
          \put(2.3,  3.2){\vector( 1, 3){ 3.1 }} 
          \put(9.5,  2.6){\vector(-1, 3){ 3.28}} 
          \put(2.0, 12.7){\vector( 0,-1){ 9.3 }} 
          \put(9.7, 12.7){\vector( 0,-1){10.1 }} 
          \put(4.8,  7.0){\vector( 4,-1){ 2.3 }} 
          \put(4.3,  6.7){\vector(-1,-2){ 1.8 }} 
          \put(7.62,5.76){\vector( 1,-2){ 1.67}} 
     \end{picture}}
\end{picture}
\caption{Preferences in the proof of Theorem~$\!$\arabic{Theorem1}.}
\label{Fig4}
\end{figure}

\begin{lem}\label{Lemma1}
In the above assumptions, {\rm(i)} $\l\_2=\l\_3$; {\rm(ii)}
$w\_2=w\_3$.
\end{lem}

\begin{pf*}{Proof of Lemma~\ref{Lemma1}.}
First, assume that $(\ww,\ll)$ is a win-loss index of $k$th order.  If
$k=0$, (i) and (ii) of Lemma~\ref{Lemma1} hold. Note also that
$w\_6=w\_7$. If $k>0$, there exists a win-loss index $(\ww',\ll')$ of
$(k-1)$th order and functions $f$, $g$, $\alpha$ and $\beta$ such that
Eqs~(\ref{18}) hold for $i=1\cdc 7$. In particular,
\medskip

$w\_6={1 \over \alpha(1,2)}
f\left(\left<(1,w'_1),(0,w'_2)\cdc (0,w'_7)\right>\right)=w\_7$.
\begin{itemize}
\item[(i)]  Similarly, $\l_2=\l_3$.
\item[(ii)] $w\_2={1 \over \alpha(1,1)}
f\left(\left<(0,w'_1)\cdc (0,w'_5),(1,w'_6),(0,w'_7)\right>\right)$,
\item[\phantom{(ii)}] $w\_3={1 \over \alpha(1,1)}
f\left(\left<(0,w'_1)\cdc (0,w'_5),(1,w'_7),(0,w'_6)\right>\right)$.
\end{itemize}

In the latter formula we interchanged the two last elements of
the multiset $W'_3$. As we have seen above, for every win-loss index of
finite order, $w\_6=w\_7$, hence $w'_6=w'_7$. Therefore, $w\_2=w\_3$.

For win-loss indices of infinite order, (i) and (ii) are proved
similarly.

Suppose now that $(\ww,\ll)$ is an arbitrary win-loss index. Then there
exists a sequence of win-loss indices of finite or infinite order
$((w^1,\l^1),(w^2,\l^2),\ldots)$ and functions $\varphi$ and $\psi$
such that
\begin{eqnarray*}
 w\_i&=&\varphi(w_i^1,\l_i^1,w_i^2,\l_i^2,\ldots),
\\
\l\_i&=&\psi(w_i^1,\l_i^1,w_i^2,\l_i^2,\ldots), \qquad i=1\cdc 7.
\end{eqnarray*}
Since for every $t=1,2,\ldots$, $\;\l_2^t=\l_3^t$ and $w_2^t=w_3^t$, we
get $\l\_2=\l\_3$ and $w\_2=w\_3$. Lemma~\ref{Lemma1} is proved.
\end{pf*}

Since $s\_2=h(w\_2,\l\_2)$ and $s\_3=h(w\_3,\l\_3)$, Lemma~\ref{Lemma1}
implies $s\_2=s\_3$.

Let us contrast alternatives 2 and 3 by {\scm}. It follows that if
either $s\_6>s\_7$ or $s\_7>s\_6$, then $s\_2=s\_3$ is impossible.
Hence, $s\_6=s\_7$. Now, contrasting alternatives 4 and 5, we get
$s\_4>s\_5$.

Finally, contrasting 6 and 7 and using the proved statements
$s\_2=s\_3$ and $s\_4>s\_5$ we get $s\_6>s\_7$, which contradicts
$s\_6=s\_7$ obtained above. Thus, our assumption is wrong, and {\scm}
is broken.
\qed
\end{pf}

Theorem~\ref{Theorem1} can be explained as follows. All functions in
Definitions~\ref{Def1}--\ref{Def4} are arbitrary, we impose no
constraints on them. The only non-arbitrariness was that ``wins" and
``losses" were {\em treated separately\/} (whereas ``draws" influenced
both win-scores and loss-scores).  This separation was not complete,
because the expressions (\ref{18}) and (\ref{19}) for win-scores
$w\_i$ also contained the total losses $c_i^{\smm}$ and vice versa;
such a cross-dependence extended to win-indices and loss-indices
(\ref{lemon}) too. In spite of that, as can be shown, a win of
alternative $j$ over alternative $k$ can affect the win-score of $i\ne
k$ only if there exists a directed path in the preference multidigraph
from $i$ to $j$ or from $i$ to $k$.  Similarly, this win of $j$ over
$k$ can affect the loss-score of $i\ne j$ only if there exists a
directed path from $j$ or from $k$ to $i$.  Therefore the ultimate
scores are sensitive only to directed paths of these two kinds. As can
be concluded from the proof of Theorem~\ref{Theorem1}, {\scm} implies
the sensitivity of the ultimate scores to the influences that spread
not only through directed paths, but also via paths with altering
directions of arrows.

\section{Win-Loss Unifying Procedures}
\label{Section8}

In this section we consider two scoring procedures that do not
separate wins and losses in the score equations, treating them {\em
uniformly}.  Another important difference from the procedures of
Section~\ref{Section6} is that here the function combining the
comparison outcomes and the estimates of the opponents (function $f$
in Eqs~(\ref{1})) has an additive (instead of productive) form. It is
worth noting that if the score of $i$ is represented through
multiplicative terms such as $a\_{ij}s\_j$, then the influence of
$s\_j$ upon the score of $i$ depends on $a\_{ij}$: the greater
$a\_{ij}$, the stronger this influence. In the extreme case where
$a\_{ij}=0$, $s\_j$ is not taken into account at all. An additive form
of scores, on the contrary, causes a uniform influence of $\sj$ upon
$\si$, regardless of $a\_{ij}$. This appears more reasonable.

\subsection{The Least-Squares Procedure}
\label{Subsect_LS}

This procedure, first proposed by Smith \cite{Smith1956} and Gulliksen
\cite{Gulliksen1956} and then investigated in \cite[etc.]{Noether1960,%
LeeflangVanPraag1971,Leake1976,KaiserSerlin1978,Vereskov1986,%
Pliner1989}, constructs a mean square approximation of the comparison
outcomes by the differences between the desired scores:
\begin{equation}
\mize_{\ss}\suml_{i,j,p}(r_{ij}^p -(\si-\sj))^2,
\label{A1}
\end{equation}
where
\begin{equation}
r_{ij}^p=a_{ij}^p-a_{ji}^p
\label{A2}
\end{equation}
is a skew-symmetric modification of the comparison outcomes; $r_{ij}^p$
are undefined whenever $a_{ij}^p$ and $a_{ji}^p$ are undefined; the sum
extends over those $i$, $j$ and $p$ for which $r_{ij}^p$ is definite.

Partial differentiation reduces the problem (\ref{A1}) to the system of
linear equations
\begin{equation}
\si={1\over m\_i}\suml_{j,p}(\sj+r_{ij}^p),\quad i=\on,
\label{A3}
\end{equation}
where $m\_i=c_i^{\smp}+c_i^{\smm}$ is the number of comparisons of $i$.

Note some parallelism of this system of equations and that
of~(\ref{16}). A multiplicative counterpart of Eqs~(\ref{A3})
is considered in \cite{Keener1993} and fits~(\ref{1.1}).

This procedure is applicable not only to indivisible preference
profiles, but to all profiles with connected preference multigraph.
Under this condition, the rank of the system~(\ref{A3}) is $n-1$, and
scores $s\_1\cdc s\_n$ can be found up to an additive constant. For
the complete preferences, this procedure reduces to that of row sums,
and it is rather well-performing for incomplete numerical (weighted)
preferences.  However, in the case of discrete incomplete preferences
which we consider here, it can produce some unnatural results (see the
proof of the following proposition).

\begin{prop}\label{Prop4}
\newcounter{Prop4}\setcounter{Prop4}{\value{thm}}
The Least squares procedure breaks {\scm}.
\end{prop}

\begin{pf}
Consider the preference multidigraph shown in Fig.~\ref{Fig5}.

\begin{figure}[htb]
\setlength{\unitlength}{5mm}
\begin{picture}(2,10)
\put(9.5,1){\begin{picture}(8,8)
     \put( 0,0){\line(1,0){8}}
     \put( 0,8){\line(1,0){8}}
     \put( 0,0){\line(0,1){8}}
     \put( 8,0){\line(0,1){8}}
          \put(1.8,6   ){$1$}
          \put(6,  6   ){$2$}
          \put(3.9,4   ){$3$}
          \put(3.9,1.74){$4$}
          \put(2.3,5.9 ){\vector( 1,-1){1.5}} 
          \put(5.9,5.9 ){\vector(-1,-1){1.5}} 
          \put(4.1,3.8 ){\vector( 0,-1){1.4}} 
          \put(2.1,5.8 ){\vector( 1,-2){1.7}} 
     \end{picture}}
\end{picture}
\caption{Preferences in the proof of Proposition~$\!$\arabic{Prop4}.}
\label{Fig5}
\end{figure}

The least-squares score of 2 is greater than that of 1, in spite of
that 1 has an extra win. This breaks {\scm}.
\qed
\end{pf}

This proof can be commented as follows. The Least-squares procedure
better fits numerical preferences; it punishes 1 for the win over 4
because it ``expects" that if 1 beats 3 and 3 beats 4 with
the same ``intensity", 1 should beat 4 with a greater intensity. The
following procedure eliminates this flaw.

\subsection{The Generalized Row Sum Procedure}
\label{Subsect_GRS}

This procedure \cite{Chebotarev1989,Chebotarev1994} can be considered
as a Bayesian modification of the previous one. It is derived
axiomatically and has Markov chain and graph theoretic interpretations
\cite{Shamis1994,ChebotarevShamis1995}. The scores satisfy the system
of equations
\begin{equation}
\si=\ve\suml_{j,p}(\gamma r_{ij}^p-(\si-\sj)),\; i=\on,
\label{A4}
\end{equation}
where $\ve$ is a positive parameter, $\ve\le(m(n-2))^{-1}$, and
$\gamma=mn+\ve^{-1}$.

\begin{prop}\label{Prop5}
The Generalized row sum procedure satisfies Self-Consistent
Monotonicity.
\end{prop}

The proof is given in \cite{ChebotarevShamis1997}. The form of these
scores suggests what kind of procedures satisfies {\scm}. A class of
such procedures is described in the following section.

\section{A Sufficient Condition of {\SCM}}
\label{Section9}

\begin{thm}\label{Theorem2}
Suppose that a scoring procedure $\s$ is such that there exists a
function $f$ defined on finite nonempty multisets of real triples and
possessing of the following properties.
\begin{itemize}
\item[{\rm(}i{\rm)}] for every $\A$,
\begin{equation}
f(\langle(a_{ij}^p,\si,\sj)\mid j,p\rangle)=0, \quad i=\on.
\label{A5}
\end{equation}
\item[{\rm(}ii{\rm)}] $f$ is
        \begin{itemize}
          \item[--] increasing in every $a_{ij}^p$,
          \item[--] increasing in every $\sj$  $(j\ne i)$,
          \item[--] decreasing in $\si$.
        \end{itemize}
\item[{\rm(}iii{\rm)}] {\em Let $V_i=\langle(a_{ij}^p,\si,\sj)\mid
j,p\rangle$.  For every $\A$ and for every $i,k\in J$},
        \begin{itemize}
          \item[] if $(1,s\_i,\sk)\in V_i$, then
          $f(V_i\setminus\left<(1,s\_i,\sk)\right>)<f(V_i)$;
          \item[] if $(0,s\_i,\sk)\in V_i$, then
          $f(V_i\setminus\left<(0,s\_i,\sk)\right>)>f(V_i)$.
        \end{itemize}
\end{itemize}
Then scoring procedure $\s$ satisfies {\scm}.
\end{thm}

\begin{pf}
Assume that the relation between $i$ and $j$ described in the statement
of {\scm} holds, but $\si<\sj$ (respectively,  $\si\le\sj$ in the
strict case). Then, by (ii) and (iii) of
Theorem~\ref{Theorem2}, the left-hand side of the $i$th equation in the
system~(\ref{A5}) exceeds the left-hand side of the $j$th equation, and
hence they both cannot be equal to zero. This contradiction proves the
theorem.
\qed
\end{pf}

In \cite{ChebotarevShamis1998} a modification of SCM is proposed which
is met by those and only those scoring procedures that have the
implicit form (i)--(ii) of Theorem~\ref{Theorem2}.

\section{Four Procedures Satisfying SCM}
\label{Section10}

Table~\ref{Table1} lists four procedures that satisfy {\scm}.

\begin{table}[htb]
\caption{Scoring procedures that satisfy SCM.}
$$
\begin{array}{lccl}
\hline
{\rm Procedure} & {\rm Domain}^* &\phantom{.}& f(\cdot)\;\mbox{of
Theorem~\ref{Theorem2}}\\
\hline
\mbox{Zermelo \cite{Zermelo1928}, Bradley\&Terry
\cite{BradleyTerry1952},...}
& {\rm ID}
&& \suml_{j,p}\left(a_{ij}^p-\frac{\textstyle\si}{\textstyle
\si+\sj}\right) \\
\mbox{Daniels \cite{Daniels1969}, Ginovker \cite{Ginovker1981}}
& {\rm ID}
&&
\suml_{j=1}^n\left(a\_{ij}\frac{\textstyle\sj}{\textstyle\si}-
a_{ji}\frac{\textstyle\si}{\textstyle\sj}\right)\\
\mbox{Cowden \cite{Cowden1975}}
& {\rm ID}
&& \suml_{j,p} (a\_{ij}\sj(1-\si)-a\_{ji}\si(1-\sj))\\
\mbox{Generalized row sum procedure %
\cite{Chebotarev1994}}
& {\rm U}
&& \ve\suml_{j,p}(\gamma r_{ij}^p-(\si-\sj))-\si\\
\hline
\end{array}
$$
$^*$ ID means indivisible preference profiles, U all profiles.
\medskip
\label{Table1}
\end{table}

The second column conveys the set of preference profiles to which the
corresponding procedure is applicable. Of the four methods, the
Generalized row sum procedure is based on the solution of a linear
system of equations, three others are reduced to nonlinear systems of
equations, which are usually solved with iterative algorithms.

A question arises of how to compare further the procedures that satisfy
{\scm}. The last section contains two remarks on the possible
additional axioms.

\section{Concluding Remarks}
\label{Section11}

In a follow-up paper we will continue the axiomatic testing of indirect
scoring procedures and turn to their axiomatic derivation. We now give
two possible conditions which can supplement {\scm}.

\subsection{Macrovertex Independence}
\label{Subsect_Macrover}

The main idea of {\scm} is sensitivity of the aggregating
procedure. However, there exists an independence condition called
{\em Macrovertex independence\/} that does not contradict SCM and seems
rather natural (cf.\ the discussion after Theorem~\ref{Theorem1}).

A subset $M$ of the set of alternatives $J$ is called a {\em
macrovertex\/} of the preference multigraph if for every $i,j\in M$ and
every $k\in J\setminus M$, $n\_{ij}=n\_{jk}$, where $n\_{ik}$ is the
number of comparisons between $i$ and $k$ in $\A$.
\medskip

{\bf Macrovertex independence}. {\em The comparison outcomes of the
alternatives within a macrovertex have no influence on the scores of
the alternatives outside the macrovertex.}

\subsection{Splitting Balance}
\label{Subsect_Splitting}

Let us note that {\scm} is not a very restrictive condition, since it
leaves room for some preconception. To conclude with what we started,
we address our first example again and give it a sport interpretation.
Suppose that Fig.~\ref{Fig1}.d depicts an incomplete chess tournament.
Assume that 4 is the world chess champion, 1, 2 and 3 being university
students.  Our aim is to estimate the strength of the players taking
into account prior information. Most probably, the world champion is
the strongest player, and we should allot the highest rank to 4.
Player 2 managed to make three draws with the champion.  That is
great, and we rank 2 second in spite of three his/her losses to 1.
Then 1 is the third and 3 the last. Note that this rank order does
not contradict SCM! This is what we meant by saying that SCM leaves
room for some preconception. There exists a logic (a bit biased) that
justifies ranking $(4,2,1,3)$. To avoid this, the following condition
may be added.
\medskip

{\bf Splitting Balance}. {\em Suppose the set of alternatives $J$ can
be split into $J_1$ and $J_2$ such that for no $i\in J_1$ and $j\in
J_2$, $a\_{ji}>0$. Then there exist $i\in J_1$ and $j\in
J_2$ such that $\si\ge\sj$.}
\bigskip

The objective of this paper was to put together various indirect
scoring procedures (mainly, of infinite order) and to subject them to
the analysis in the spirit of the social choice theory. We considered
the case of incomplete preferences and an axiom we referred to as
{\scm}. It was established that many of the known indirect scoring
procedures are {\em win-loss combining procedures\/} and, by
Theorem~\ref{Theorem1} of Section~\ref{Section7}, they break {\scm}. In
our opinion, some of them are applicable in the analysis of social
networks, since here {\scm} is not so attractive as in preference
aggregation (cf. \cite{Katz1953,Taylor1969,Friedkin1991}). For
instance, the popularity index of an individual should not strongly
depend on his/her own responses.

In Theorem~\ref{Theorem2} we proposed a sufficient condition for {\scm}
and then listed four procedures that satisfy it.  The analysis of
indirect scoring procedures started in this paper is intended to be
continued with additional axioms (two possible conditions were
mentioned in this section), and with the end goal of axiomatic
construction of such procedures.


\end{document}